\documentclass[reqno]{amsart}

\usepackage{a4wide}

\usepackage{times}
\usepackage[T1]{fontenc}
\usepackage{pstricks}
\usepackage{amssymb}

\newtheorem{lem}{Lemma}
\newtheorem{thm}{Theorem}
\newtheorem{cor}{Corollary}

\renewcommand{\ge}{\geqslant}
\renewcommand{\le}{\leqslant}

\newcommand{\al}{\alpha}
\newcommand{\e}{\varepsilon}

\newcommand{\R}{{\mathbb{R}}}

\newcommand{\N}{{\mathbb{N}}}
\newcommand{\T}{{\mathbb{T}}}
\newcommand{\Z}{{\mathbb{Z}}}

\newcommand{\0}{{\mathbf{0}}}

\newcommand{\q}{{\mathbf{q}}}

\newcommand{\rb}{\mathbf{r}}

\newcommand{\col}{\colon}

\newcommand{\DA}{Diophantine approximation}

\newcommand{\HD}{Hausdorff dimension}
\newcommand{\HM}{Hausdorff measure}

\newcommand{\LM}{Lebesgue measure}

\newcommand{\dist}{\operatorname{dist}}
\newcommand{\codim}{\operatorname{codim}}
\newcommand{\Ja}{Jarn\'{\i}k}
\newcommand{\JB}{Jarn\'{\i}k-Besicovitch}

\newcommand{\cEk}{{\chi}_{_{E_k}}}
\newcommand{\cEkl}{{\chi}_{_{E_k\cap E_l}}}
\newcommand{\cEl}{{\chi}_{_{E_l}}}
\newcommand{\cEn}{{\chi}_{_{E_n}}}

\newcommand{\ub}{\mathbf{u}}

\newcommand{\xb}{\mathbf{x}}
\newcommand{\eg}{{\it{e.g.}}}

\newcommand{\ip}{\,{\cdot}\,}
\newcommand{\ie}{{\it{i.e.}}}
\newcommand{\im}{infinitely many}

\newcommand{\cC}{\mathcal{C}}
\newcommand{\cH}{\mathcal{H}}
\newcommand{\cP}{\mathcal{P}}

\newcommand{\vphi}{\varphi}
\DeclareMathOperator{\Card}{Card}

\newcommand{\cN}{\mathcal{N}}
\newcommand{\cR}{\mathcal{R}}

\newcommand{\tB}{\widetilde{B}}

\newcommand{\KG}{Khintchine-Groshev} 
\newcommand{\K}{Khintchine}

\newcommand{\Spr}{Sprind{\v z}uk}

\def\tB{\widetilde{B}}
 
 \def\Bqr{B(q;\rho)}
\def\Bqrp{B(q';\rho')}
 \def\Bqp{B(q;\psi(q))}
\def\Bqbd{B(\q;\delta)}
\def\Bqbr{B(\q;\rho)}
\def\Bqbrp{B(\q';\rho')}
\def\B|qb|p{B(\q;\psi(|\q|))}

\def\Rq{R_{\q}}
\def\tp{\widetilde{\psi}}
\def\tr{\tilde{\rho}}

\title{Diophantine approximation, Khintchine's theorem, torus geometry and
                       Hausdorff  dimension}

\author{M. M. Dodson}
\address{
  University of York\\
  Heslington \\
  York \\
  YO10 5DD \\
  UK
}
\email{mmd1@york.ac.uk}

\begin{document}

\begin{abstract}
  A general form of the Borel-Cantelli Lemma and its connection with
  the proof of Khintchine's Theorem on Diophantine approximation and
  the more general \KG{} theorem are discussed. The torus geometry in
  the planar case allows a relatively direct proof of the planar
  Groshev theorem for the set of $\psi$-approximable points in the
  plane.  The construction and use of Haudsorff measure and dimension
  are explained and the notion of ubiquity, which is effective in
  estimating the lower bound of the \HD{} for quite general lim sup
  sets, is described. An application is made to obtain the \HD{} of
  the set of $\psi$-approximable points in the plane when
  $\psi(q)=q^{-v}$, $v>0$, corresponding to the planar \JB{} theorem. 
\end{abstract}

\maketitle

\section{Diophantine approximation}

\DA{} is a quantitative analyis of the density of the rationals in the
reals.  It is easy to see from the distribution of the integers $\Z$
in the real line $\R$, that given any $\al\in\R$ and any $q\in\N$,
there exists a $p=p(\al,q)\in\Z$ such that
\begin{equation*}
  |q\al-p|\le 1/2 \ {\text{or} } \ |\al-p/q|\le \frac1{2q}.
\end{equation*}
It is possible to do better using continued
fractions (see~\cite{Casselshort,HW}) and thanks to Dirichlet's box
argument~\cite{HarmanMNT}, to obtain a best possible result.
\begin{thm}[Dirichlet]
 Given $\alpha\in\R$, $N\in\N$, there exist integers 
$p, q$ with $1\le q\le N$ such that
\begin{equation*}
   |\al-p/q|\le \frac1{q(N+1)}.
\end{equation*}
  \end{thm}
Given any $\xi\in\R$, the  convenient notation
\begin{equation*}
\|\xi\|:=\min\{|\xi - p|\colon p\in \Z\}   
\end{equation*}
will be used.  Thus Dirichlet's theorem implies that given any $\al\in\R$,
there are infinitely many $q\in\N$ such that
\begin{equation*}
  \|q\al\|= \min\{|q\alpha-p|\colon p\in\Z\} <\frac1{q}.
\end{equation*}
More generally, an error term or {\em approximation} function $
\psi\colon\N\to (0,\infty)$, where $\lim_{q\to\infty} \psi(q)=0$, is
introduced and the solubility of
\begin{equation}
\label{eq:psineq}
   \|q\al\| <  \psi(q)
\end{equation}
considered ($\psi(q)=1/q$ in Dirichlet's theorem).  Note that although
restricting the approximation to rationals $p/q$ with $(p,q)=1$ is
natural and indeed is associated with the Duffin-Schaeffer conjecture
(see~\cite{DuffinSchaeffer}), coprimality does not arise in the present
formulation.

The point $\alpha$ is said to be {\bf{{\em $\psi$-approximable}}}
if~\eqref{eq:psineq} holds for infinitely many $q\in\N$.  The set
$W(\psi)$ of $\psi$-approximable numbers is invariant under
translation by integers and so there is no loss of generality in restricting attention to the unit interval and considering
   \begin{align*}
     W(\psi):&=\{\al\in[0,1]\colon \|q\al\|<\psi(q)
     \text{ for infinitely many } q\in\N\}\\
&=\{\al\in[0,1]\colon |\al-p/q|<\psi(q)/q
     \text{ for infinitely many } p\in\Z, q\in\N\}\\
     & = \bigcap_{N=1}^\infty\bigcup_{q=N}^\infty\bigcup_{p=0}^{q}
     \left(\frac{p}{q}- \frac{\psi(q)}{q}, \frac{p}{q}+\frac{\psi(q)}{q}\right)
\cap [0,1].
   \end{align*}
   The set $W(\psi)$ is a lim sup set, as can be seen by letting for
   each $q\in\N$ and $\rho>0$,
\begin{equation*}
\Bqr =\bigcup_{p=0}^{q} \left(\frac{p}{q}- 
\frac{\rho}{q},\frac{p}{q}+\frac{\rho}{q}\right)
\cap [0,1],
\end{equation*}
so that $\Bqr$ is a $\rho$-neighbourhood of the resonant set (so
called from the connection with the physical phenomenon of resonance)
   \begin{equation*}
     R_q:=\left\{0,\frac1{q}, \dots, \frac{p}q{}, \dots, 
\frac{q-1}{q}, 1 \right\}
   \end{equation*}
   and
   \begin{align}
\label{eq:Wp}
         W(\psi) &= \bigcap_{N=1}^\infty\bigcup_{q=1}^{N}
     \Bqp=\limsup_{q\to\infty} \Bqp.
   \end{align}
Moreover, $|\Bqr|=2\rho$, since the set of points in $[0,1]$
   satisfying
\begin{equation*}
  \left|\al-\frac{p}{q}\right| < \frac{\rho}{q},
\end{equation*}
shown in Figure~1, has length
\begin{equation*}
  \frac{\rho}{q} + (q-1) \frac{2\rho}{q} +   \frac{\rho}{q} 
=   2\rho.
\end{equation*}
  \begin{center}
  \begin{pspicture}(-4,0.5)(4,2)
\psline(-4,1.5)(4,1.5) 
  \rput(-4,1.25){$0$}
  \rput(4,1.25){$1$}
  \rput(-2,1.2){$\frac{1}{q}$}
  \rput(2,1.2){$\frac{q-1}{q}$}

\psdots[dotstyle=|](-4,1.5)(-2,1.5)(0,1.5)(2,1.5)(4,1.5)
\psline(-4,1.52)(-3.5,1.52) 
\rput(-3.5,1.2){$\frac{\rho}{q}$}
\psline(-2.5,1.52)(-1.5,1.52)
\psline(-.5,1.52)(.5,1.52)
\psline(1.5,1.52)(2.5,1.52)
\psline(3.5,1.52)(4,1.52)
\end{pspicture}
  \end{center}

  \begin{center}
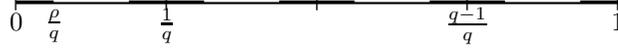
\small{Figure 1. Points $p/q$ and neighbourhood
  $\Bqr$ in $[0,1]$.}
\end{center}
\vspace{0,1in} We are interested in the `size' of $W(\psi)$.  The size
of a measurable set $E$ will be interpreted as its Lebesgue measure,
denoted by $|E|$.  The question of the measure of $W(\psi)$ is almost
completely answered by
\subsection{Khintchine's theorem}
\begin{thm}
  The Lebesgue measure of $W(\psi)$ is given by
  \begin{equation*}
    |W(\psi)| = \begin{cases} 0, 
           &\text{if } \sum_{k=1}^\infty \psi(k)<\infty, \\
1,
          &\text{if }  \sum_{k=1}^\infty \psi(k)=\infty 
\ \text{and}\ \psi(k) \ \text{is non-increasing}.
                  \end{cases}
  \end{equation*}
\end{thm}
The theorem corresponds to our intuition, as when the approximation
function $\psi$ is large,  there is a better chance of the
inequality being satisfied and vice-versa (for more details,
see~\cite{Casselshort,Sprindzuk}).  Thus the Lebesgue measure of
$W(\psi)$ is $1$ when $\psi(q) = 1/(q \log q)$ and $0$ when $\psi(q) =
1/(q(\log q)^{1+\e})$ for any positive $\e$.  This `0-1' property is a
feature of the metrical theory and reflects its links with
probability.  Indeed the result is reminiscent of the Borel-Cantelli
Lemma from probability theory (see~\cite{KandT}).  Let $E_j, j=1,2
\dots$ be a sequence of events in a probability space $(\Omega, P)$,
let
\begin{align*}
E&=\{x\in\Omega\colon x\in E_j \text{ infinitely often }\} \\
&=\bigcup_{N=1}^\infty \bigcap_{r=N}^\infty
  E_r=\limsup_{N\to\infty} E_N.
\end{align*}
Then the lim sup $E$ is the the set of points lying in infinitely many
events $E_j$ and $P(E)$ is the probability that infinitely many events
$E_j$ occur. 
\begin{lem}[Borel-Cantelli]
  \label{lem:BC}
  \begin{equation*}
P(E)= \begin{cases} 0 , &\text{if } \sum_{j=1}^\infty P(E_j)<\infty, \\
                 1, &\text{if } \ E_j \ \text {totally independent and
                 } \ \sum_{j=1}^\infty P(E_j)=\infty.
                 \end{cases}
  \end{equation*}
\end{lem}
Borel proved the lemma assuming total independence of the events, when 
for any distinct events $E_{j_1},\dots,E_{j_k}$,
\begin{equation*}
P (E_{j_1}\cap  \dots
\cap E_{j_k})=P(E_{j_1}) \dots  P(E_{j_k}).
\end{equation*}
Cantelli observed that independence is not needed when the sum of
probabilities converges.  In the case of divergence, the result holds
under total {\it quasi}-independence, when for some $K\ge1$,
\begin{equation*}
P (E_{j_1}\cap  \dots
\cap E_{j_k})\le K P(E_{j_1}) \dots  P(E_{j_k}).
\end{equation*}

If Lebesgue measure is interpreted as probability, then
$|\Bqp|$ 
corresponds to the probability that a point $\al\in [0,1]$ falls into
$\Bqp$ and the proof of \K's theorem in the case of convergence is
essentially that in the Borel-Cantelli Lemma. Indeed, the set of
points in $[0,1]$ satisfying~\eqref{eq:psineq} for a given $q\in\N$
and $\rho=\psi(q)$, shown in Figure~1, has length
\begin{equation*}
  \frac{\psi(q)}{q} + (q-1) \frac{2\psi(q)}{q} + \frac{\psi(q)}{q} 
=   2\psi(q). 
\end{equation*}
Hence for any $N\in \N$, 
\begin{equation*}
  |W(\psi)|\le 2 \sum_{q=N}^\infty  \psi(q) \to \ 0 \ {\text{as}} \ N\to\infty,
\end{equation*}
whence $|W(\psi)|=0$. Sets of Lebesgue measure 0 will be called
\emph{null}. 

The case of divergence is more difficult.  The Borel-Cantelli Lemma
assumes total independence to deal with the divergence case and so is
useless.  However a more general lower bound for lim sup set of the
sets $E_j$ is available. It suits our purposes to express the result
in terms of the Lebesgue measure $|E_j|$ of the sets $E_j$.
\begin{thm}  
\label{thm:pwi2}
Let $E_j$, $j=1,2,\dots$, be a sequence of Lebesgue measurable sets in
$\Omega=[0,1]^n$  
and suppose that 
 \begin{equation}
\sum_{j=1}^\infty |E_j| =\infty.
   \label{eq:Esum}
 \end{equation}
 Then the Lebesgue measure of $E:=\limsup_{N\to\infty}
 E_N:=\bigcap_{N=1}^\infty \bigcup_{k=N}^\infty E_k$ satisfies
\begin{equation}
 \label{eq:muE}
 |E|  \ge  \limsup_{N\to\infty}\frac{(\sum_{k=1}^N |E_k|)^2}
 {\sum_{k=1}^N\sum_{l=1}^N |E_k\cap E_l|}.
\end{equation}
\end{thm}
Because of the importance of the result, a proof will be given.  This
is based on `mean and variance' arguments (see~\cite{Sullivan82} and
also~\cite{HarmanMNT,Sprindzuk}).  It appears as an exercise
in~\cite{Chung} and there are numerous variants, \eg,~\cite{KS64}.
\begin{proof}
  For each $n=1,2,\dots$, let $\nu_n\colon \Omega\to [0,\infty]$ be
  the counting function of the number of $E_j$ into which $x$ falls,
  so that
$$
\nu_n(x):=\sum_{k=1}^n \cEk(x)\le \lim_{n\to\infty}\nu_n(x) 
=\sum_{k=1}^\infty \cEk(x):=\nu(x). 
$$
 The average of $\nu_n$ over $\Omega$ is given by 
  \begin{equation}
    \label{eq:an}
 A_n =\int_\Omega \nu_n(x)\, dx = \sum_{k=1}^n \int_\Omega \cEk(x)\,dx =
  \sum_{k=1}^n |E_k|
  \end{equation}
  and 
 \begin{equation}
    \label{eq:int1}
   \int_\Omega \frac{\nu_n(x)}{A_n}\,dx
=\frac{\int_\Omega \nu_n(x)\,dx}{ \int_\Omega \nu_n(x)\, dx}=1.
 \end{equation}
 Thus the average of $\nu_n(x)/A_n:=f_n(x)$ over $\Omega$ is 1.
 By~\eqref{eq:Esum} and~\eqref{eq:an}, $A_n\to\infty$ as $n\to\infty$
 and so $\nu(x)$ can be infinite.  Now  the lim sup $E$ is given by
\begin{align*}
  E & = 
   \{x\in \Omega \col \sum_{n=1}^\infty \cEn(x)=\infty\} = \{x\in
  \Omega \col \nu(x)=\infty\}, 
\end{align*} 
so that 
 \begin{equation*}
\nu(x)=\sum_{k=1}^\infty \cEk(x)<\infty \ \text{if and only if} \ x\notin E,
 \end{equation*}
 \ie, $\nu(x)<\infty$ if and only if $x\in \Omega \setminus E= E^c$.  But by
 hypothesis, $A_n\to\infty$ as $n\to\infty$, whence for each $x\in
 E^c$, $f_n(x)=\nu(x)/A_n\to 0$ as $n\to\infty$.  Thus the
 contribution to $\int_\Omega f_n=\int_E f_n+\int_{E^c} f_n=1$ is
 mainly from $E$.  Assume for the moment that, as one would expect,
 \begin{equation}
\label{eq:limintfn}
\int_E f_n =\frac1{A_n} \int_E\nu_n \to 1 \ \text{as} \ n\to\infty, 
 \end{equation}
so that from Cauchy's inequality, 
\begin{align*}
1   &\le \frac{\left(\int_E \nu_n(x)^2 \,dx \right)^{1/2}|E|^{1/2}}{A_n}
  + o(1),
\end{align*}
whence on rearranging,
\begin{equation*}
\limsup_{n\to\infty}  \frac{A_n^2}{\int_E \nu_n(x)^2 \,dx}\ge 1. 
\end{equation*}
But 
\begin{align*}
\int_\Omega \nu(x)^2 \,dx & =   
\int_\Omega \sum_{k=1}^n \cEk(x)
\sum_{l=1}^n \cEl(x)\, dx = 
\int_\Omega  \sum_{k,l=1}^n \cEk(x) \cEl(x)\,dx \\ 
                    &= \int_\Omega  \sum_{k,l=1}^n \cEkl(x) \, dx 
=   \sum_{k,l=1}^n |E_k\cap E_{l}|,
\end{align*} 
and the result follows.

To prove that $\int_{E^c} f_n=o(1)$
and hence prove~\eqref{eq:limintfn}, consider the sequence
$(f_n=\nu_n/A_n\colon E^c\to \R)$ of (measurable) non-negative
functions.
This sequence is well defined with $\lim_{n\to\infty} f_n(x)=0$ for
each $x\in E^c$.  Further $|E^c|\le1$.  Hence by Egeroff's theorem,
given $\eta>0$, there exists a measurable subset $F_\eta\subseteq E^c$
such that $|F_\eta|<\eta$ and $f_n\to 0$ uniformly on $E^c\setminus
F_\eta$.  Hence
\begin{align*}
  \lim_{n\to\infty} \int_{E^c} f_n(x) \, dx&=  
  \lim_{n\to\infty} \int_{E^c\setminus F_\eta} f_n(x) \, dx + 
  \lim_{n\to\infty} \int_{ F_\eta} f_n(x) \, dx \\
 &=   \int_{E^c\setminus F_\eta}\lim_{n\to\infty} f_n(x) \, dx + 
  \lim_{n\to\infty} \int_{ F_\eta} f_n(x) \, dx\\
 &= 0 +\lim_{n\to\infty} \int_{ F_\eta} f_n(x) \, dx.
\end{align*}
But $\int_\Omega f_n=1$ for each $n$, whence $\int_{F_\eta} f_n\to 0$ as
$\eta \to 0$. Hence
$$
\lim_{n\to\infty} \int_{E^c} f_n(x) dx =\lim_{n\to\infty} \int_{E^c} \frac{\nu_n(x)}{a_n} dx = 0, 
$$
as required
\end{proof}

\subsection*{Pairwise independence}
The sets $E_k$, $E_l$ are {\em pairwise quasi independent} if there
exists a constant $C$ such that for all distinct $k,l$,
 \begin{equation*}
   |E_k\cap E_l| \le C |E_k||E_l| 
\end{equation*}
and are {\em pairwise independent} if for all distinct $k,l$,
 \begin{equation*}
   |E_k\cap E_l| = |E_k||E_l|. 
\end{equation*}
The results below are immediate consequences of~\eqref{eq:muE}
and~\eqref{eq:Esum}. 
 \begin{cor}
\label{cor:pwqi}
   If the $E_j$ are pairwise quasi independent, then $|E|>0$.
\end{cor}
 \begin{cor}
\label{cor:pwi}
   If the $E_j$ are pairwise  independent, $E=1$.
\end{cor}

It can be shown that the $\Bqr$
are pairwise quasi-independent, \ie, for distinct $q,q'$,
\begin{equation*}
  |\Bqr\cap \Bqrp|\le K |\Bqr| |\Bqrp|,
\end{equation*}
so that by Corollary~\ref{cor:pwqi}, $|W(\psi(q)|>0$.  Establishing
pairwise quasi independence involves some lengthy and difficult
technicalities and will be omitted (proofs are given
in~\cite{Sprindzuk,Sullivan82}).  However, once established, the full
result follows from an ergodic-type theorem of
Gallagher~\cite{Gallagher65} or from Lebesgue density that $|W(\psi)|=1$.
This `all or nothing' or `$0$-$1$' law was originally proved by \K{}
(who used continued fractions~\cite{KhCF}, limiting the proof to
$\R$).  Other proofs based on pairwise quasi-independence or mean and
variance arguments (see for example~\cite[Chaper~VII]{Casselshort}) in
conjunction with density or ergodic ideas can be extended to higher
dimensions and other generalisations. Recently Cornelia Drutu has used
pairwise quasi-independence for \DA{} in a symmetric spaces
setting~\cite{drutu05}.

\section{Higher dimensions}
In higher dimensions ($\R^n$), there are two natural forms of
Diophantine approximation. First, in the simultaneous form, one
considers the set of points $\xb =(x_1,\dots,x_n)\in \R^n$ for which
the inequality 
\begin{equation}
\label{eq:qxb}
\|q\xb\|:=\max\left\{||qx_1||,\dots,||qx_n|| \right\} < \psi(q)  
\end{equation}
holds for infinitely many positive integers $q$. There is also the {\it dual}
form in which one considers the proximity of the point
$\xb=(x_1,\dots,x_n)\in\R^n$ to the hyperplane 
\begin{equation}
  \label{eq:Rpq}
\{\ub\in\R^n\colon \q\cdot\ub=p\},
\end{equation}
where $p\in\Z, \q\in\Z^n$. More precisely, one considers the
solubility of the inequality
\begin{equation}
  \label{eq:2}
|\q\cdot \xb-p|<\psi(|\q|),
\end{equation}
where for each $\ub\in\R^n$, $|\ub|:=|\ub|_\infty =
\max\{|u_1|,\dots,|u_n|\}$, for \im{} $p, \q$.  This notation should
not be confused with that for the Lebesgue measure of a set. For
convenience, introduce the resonant set
\begin{equation}
\label{eq:Rqp}
   R_{\q}(p):=\{\ub\in [0,1]^n\colon |\q\ip \ub-p| = 0\}, \ p\in \Z
\end{equation} 
and denote the collection of resonant sets associated with $\q$ by
\begin{equation}
\label{eq:Rq}
  R_{\q}:=\{\ub\in [0,1]^n\colon \|\q\ip \ub\| = 0\} 
=\bigcup_{p} R_{\q}(p)
\end{equation} 
and let $\Bqbr$ be the $\rho$-neighbourhood 
$\{\ub\in [0,1]^n\colon \|\q\ip \ub\| < \rho\}$ of $R_\q$.
The above two forms of approximation are special cases of a system of
linear forms, discussed in~\S\ref{subsec:KG}.  Analytical ideas,
including Fourier series, play an important part in metrical
Diophantine approximation, as do geometrical ones, particularly so in
the dual form of \DA.

\subsection{The \K-Groshev Theorem}
\label{subsec:KG}
\K's theorem has a very general extension, originally proved by
A.~V.~Groshev~\cite{Sprindzuk} with a stronger monotonicity condition,
which includes as special cases simultaneous Diophantine approximation
and its dual, as mentioned above.  It treats real $m \times n$
matrices $X=(x_{ij})$, regarded as points in $\R^{mn}$, which are
$\psi$-{\it{approximable}}, \ie, which satisfy
\begin{equation}
  \label{eq:3}
\| \q X \| < \psi(|\q|),   
\end{equation}
for infinitely many $\q \in \Z^m$, where $\q X$ is a vector of the
following linear forms 
$$
(q_1x_{11} + \dots + q_m x_{m1}, \dots,q_1
x_{n1} + \dots + q_m x_{mn})
$$ 
and $\| \ub \|= |(\|u_1\|,\dots,\|u_n\|)|=
\max\{\|u_1\|,\dots,\|u_n\|\}$.  As the set of $\psi$-approximable
points is translation invariant under integer vectors, we can restrict
attention to the $mn$- dimensional torus $\T^{mn}$, \ie., the
$mn$-dimensional unit cube with opposite sides identified. The set of
$\psi$-approximable points in $\T^{mn}$ will be denoted by
\begin{equation*}
W(\psi;m,n) = \{X\in \T^{mn}:\|\q X \| < \psi(|\q|) {\text { for
    infinitely many }} \q \in \Z^m \}.  
\end{equation*}
To avoid complicated notation, the dependence of $W(\psi;m,n)$ on $m,n$
will usually be omitted.  
\begin{thm} 
  The $mn$-dimensional Lebesgue measure of $W(\psi)$ is given by
 \begin{equation*}
    |W(\psi)| = \begin{cases} 0, 
      &\text{if } \sum_{k=1}^\infty k^{m-1}\psi(k)^n<\infty, \\
      1,
      &\text{if }  \sum_{k=1}^\infty k^{m-1}\psi(k)^n=\infty 
      \ \text{and when } \ m=1,2 \ \ \psi(k) \ \text{is decreasing}.
                  \end{cases}
  \end{equation*}
\end{thm}
The proof is straightforward when the `volume' sum 
\begin{equation}
  \label{eq:volsum}
\sum_{k=1}^\infty  k^{m-1}\psi(k)^n
\end{equation}
converges, as the fact that $W(\psi)$ can be expressed as a lim sup
set again provides a direct and simple proof that $W(\psi)$ has
measure $0$.  However not surprisingly, as in one dimension, the
case of divergence is much more difficult and the more general
lower bound for lim sup sets is used.

\subsection{Torus geometry in the plane}

Geometrical ideas play a particularly important role in the dual form
of \DA, in which the proximity of the point
$\xb=(x_1,\dots,x_m)\in\R^m$ to the hyperplane
\begin{equation}
  \label{eq:1}
 \{\ub\in\R^m\colon \q\cdot\ub=p\},
\end{equation}
where $p\in\Z, \q\in\Z^m$ is considered. More precisely, one considers the
solubility of the inequality
\begin{equation}
  \label{eq:2a} 
|\q\cdot \xb-p|<\psi(|\q|),
\end{equation}
where $|\q|=\max\{|q_j|\colon j=1,\dots,m\}$, for \im{} $ \q$.  
 
\vspace{0.1in}
\vspace{0.1in}

    \begin{center}
  \begin{pspicture}(-2,-1)(2,3.5) 
\psset{unit=0.6cm}

\psline(-3.25,0)(4,0) 
 \rput(-2.5,-0.4){$(0,0)$}
 \psline(-2.5,-.2)(-2.5,5.7) 

\psline(-2.5,-.8)(-2.5,-1.4) 

\psline[linestyle=dashed](-3.5,5)(-1.5,6)  
\psline(-3.5,4.5)(-.75,6)  
\psline(-3.5,4.5)(-.75,6)  
\psline[linestyle=dashed](-3.5,4)(0,6)  

\psline[linestyle=dashed](-3.5,2.5)(3.5,6) 
 \psline(-3.5,2)(4.5,6) 
 \psline(-3.5,2)(4.5,6) 
 \psline[linestyle=dashed](-3.5,1.5)(4.5,5.5) 

\psline[linestyle=dashed](-3.5,0)(4.5,4) 
 \psline(-3.5,-.5)(-3.3,-.4) 
 \psline(-2.6,-0.05)(4.5,3.5) 
 \psline(-3.5,-.5)(-3.3,-.4) 
 \psline(-2.6,-0.05)(4.5,3.5) 
 \psline[linestyle=dashed](-3.55,-1.025)(-3.0,-.75) 
 \psline[linestyle=dashed](-2,-0.25)(4.5,3) 

\psline[linestyle=dashed](-.5,-1)(4.5,1.5)  
\psline(.5,-1)(4.5,1)  
\psline(.5,-1)(4.5,1)  
\psline[linestyle=dashed](1.5,-1)(4.5,.5)  
\end{pspicture}
\end{center}
\begin{center}
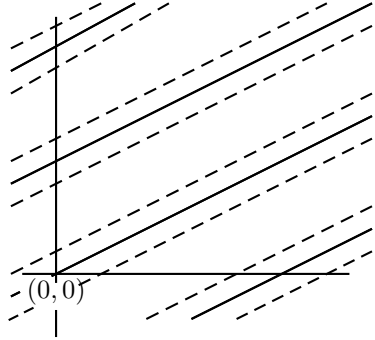
  \small{Figure 2. Resonant sets $R_{(1,-2)}$ (bold lines) and
    boundaries (dashed lines) \\
of the strips $B_{\psi(2)}(1,-2)$ in
    $\mathbb{R}^2$}
\end{center}

It turns out that there is a precise correspondence between
probabilistic independence and (algebraic) linear independence.  To
illustrate these ideas, the \KG{} theorem will be 
considered in detail
 for the planar case, where they are particularly clear.
\begin{thm} 
  The Lebesgue measure of 
  \begin{equation*}
    W(\psi):=\{\ub\in [0,1]^2\colon \|\q\ip\ub\|<\psi(|\q|) 
\text{ for infinitely many } \q\in\Z^2\}
  \end{equation*}
 satisfies
  \begin{equation*}
    |W(\psi)| = \begin{cases} 0, 
           &\text{if } \sum_{k=1}^\infty k\psi(k)<\infty, \\
1,
          &\text{if }  \sum_{k=1}^\infty k\psi(k)=\infty 
\ \text{and}\ \psi(k) \ \text{non-increasing}.
                  \end{cases}
  \end{equation*}
\end{thm}

Suppose $\q \ne \mathbf{0}$, with say $q_1 \ne 0$.  Then the resonant
set $R(\q)$ is a set of $q_1$ parallel lines in $\T^2$, a distance
$1/|q_1|$ apart in the $x_1$ direction (Figure~2).  These define
$|q_1|$ strips $S$ (from the top of a shaded strip to the top of the
adjacent shaded strip) and $\T^2 = \bigcup S$.  The set
$B(\q,\rho)$ of shaded strips $ \tilde S$ each of length $\rho/q_1$
(in the $x_1$ direction), and the ratio
\begin{equation*}
|B(\q,\rho))|:|\T^2|=|\bigcup \tilde S|:|\bigcup S| 
= |\tilde S|:|S| = 2\rho : 1
\end{equation*}
(see Figure 3).  
\begin{center}
\begin{pspicture}(-7,-1)(6,6) 
\psset{unit=0.5cm}
\psline(-5,10)(5,10) 
\psline(-5,0)(5,0) 
 \psline(-5,0)(-5,10) 
 \rput(-5,-0.5){$(0,0)$}
 \rput(5,-0.5){$(1,0)$}
 \rput(-5,10.5){$(0,1)$}
 \rput(5,10.5){$(1,1)$}
 \psline(5,0)(5,10) 
\rput(0,-.5){$x_1$} 
\rput(-5.75,5.7){$x_2$} 

\psline(-5,9)(-3,10)  
\pspolygon[fillstyle=vlines](-5,9)(-5,10)(-3,10)

\psline(-5,6)(3,10) 
 \psline(-5,5)(5,10) 
\rput(-.75,5){$S$}

\pspolygon[fillstyle=vlines](-5,4)(5,9)(5,10)(3,10)(-5,6)

\psline(-5,1)(5,6) 
 \psline(-5,0)(5,5) 
 \psline(-3,0)(5,4) 
\pspolygon[fillstyle=vlines](-5,0)(-5,1)(5,6)(5,4)(-3,0)

\psline(-5,9)(-3,10)  

\pspolygon[fillstyle=vlines](-5,9)(-5,10)(-3,10)

\psline(-5,6)(3,10) 
 \psline(-5,5)(5,10) 
 \psline(-5,4)(5,9) 

\pspolygon[fillstyle=vlines](-5,4)(5,9)(5,10)(3,10)(-5,6)

\psline(-5,1)(5,6) 
 \psline(-5,0)(5,5) 
 \psline(-3,0)(5,4) 

\pspolygon[fillstyle=vlines](-5,0)(-5,1)(5,6)(5,4)(-3,0)

\psline(3,0)(5,1)  
\pspolygon[fillstyle=vlines](3,0)(5,1)(5,0)

\end{pspicture}
 \end{center}
 \begin{center}
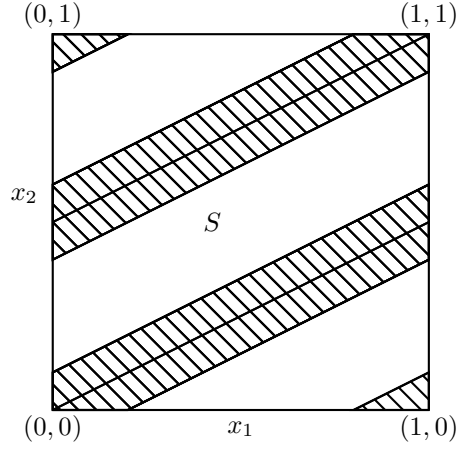
\small{Figure 3.  Strips $B_{\rho}(1,-2)$ in $\mathbb{T}^2$}
\end{center}
Thus
\begin{equation}
  \label{eq:|B|}
|\Bqbr|
=2\rho = |(-\rho,\rho)|,
\end{equation}
as in the case $m=1$.  The extension to the general case 
\begin{equation*}
  \label{eq:|Bmn|}
|\Bqbr| 
=2^n\rho^n = |(-\rho,\rho)|^n,
\end{equation*}
follows by considering $n $ copies of the 2-dimensional space spanned
by $\q$ and $\q'$ and the volume of the corresponding $mn$-dimensional
prisms.  The determination of the Lebesgue measure of $W(\psi)$ in the
case of convergence follows readily.  Fourier series can also be used
(see~\cite{Sprindzuk}).

\section{Convergence and measure 0}
The `probabilistic' interpretation discussed above, in which the
Lebesgue measure $|B(\q,\psi(|\q|))|$ of the set $B(\q,\psi(|\q|))$ is
interpreted as the probability that $\xb \in \T^n$ satisfies $\|\q\ip
\xb \|<\psi(|\q|)$, reduces the result to the Borel-Cantelli Lemma.  
For convenience, write 
\begin{equation}
  \label{eq:Bq}
  B_\q:= \B|qb|p 
\end{equation}
and consider the sets $B_\q$, $\q\in\Z^2\setminus\{0\}$, as a
 sequence $E_r$, $r=1,2,\dots$ in $\T^2$ by ordering the integer vectors
 $\q$, so that $\q=\q(r)$ is the $r$-th vector in
 $\Z^m\setminus\{\0\}$. 
 Then by~\eqref{eq:|B|} 
 \begin{equation*}
 |B_\q| =  2\psi(|\q|) ,
 \end{equation*}
so that
\begin{align*}
  \sum_{r=1}^\infty |E_r| &= \sum_{r=1}^\infty |B_{\q(r)}|
  =
  \sum_{r=1}^\infty 2\psi(|\q(r)|) 
 = 2\sum_{k=1}^\infty\sum_{|\q|=k}\psi(|\q|)=
  \sum_{k=1}^\infty\psi(k)\sum_{|\q|=k}1 \\
  &\asymp \sum_{k=1}^\infty k\psi(k),
\end{align*}
since there are $2(2k+1)\asymp k^2$ non-zero integer vectors with
$|\q|=k$ (positive quantities $a,b$ are comparable, denoted by
$a\asymp b$, if there are constants $K,K'$ such that $a\le Kb$ and
$b\le K'a$).  The convergence of the volume sum~\eqref{eq:volsum} thus
implies the convergence of the measure sum $ \sum_{r=1}^\infty |E_r|$
and hence that $|E|=|W(\psi)|=0$.  It is clear that the proof extends
to the general case.

\section{Divergence and full measure}

In the harder case when \eqref{eq:volsum} diverges, it turns out that
when $m\ge2$, the pairwise `probabilistic' independence of sets is
associated with linearly independent pairs of integer vectors, \ie,
pairs of vectors which are not collinear with the origin. Thus the
more general version of the divergence part of the Borel-Cantelli
Lemma (Corollary~1 to Theorem~\ref{thm:pwi2}) can be used.  
In this argument, the monotonicity condition can be relaxed when $m\ge
3$.
Gallagher~\cite{Gallagher65} has shown that, under a weak coprimality
condition, monotonicity can be dropped when $n\ge 2$ and $m=1$ (the
simultaneous case). Indeed even more general results, where the
argument of the error function is the vector $\q$ rather than its
supnorm $|\q|$, for $m+n>2$ were obtained for primitive solutions $\q$
by \Spr~\cite{Sprindzuk} and Schmidt~\cite{HarmanMNT}. Note that the
Duffin-Schaeffer conjecture holds for the simultaneous
case~\cite{PV90}.
Consider another vector $\q'$ and suppose that $\q, \q' \in \Z^2$ are
linearly independent (see Figure~4).  Then the set $B(\q,\rho)\cap
B(\q',\rho')$ tessellates $\T^2$ into $|\q\times \q'|$ parallelograms
$\Pi$ say, each of area $1/|\q \times \q'|$. Thus
$$
\T^2 = \cup \Pi.
$$
In addition, the set $B(\q;\rho) \cap B(\q',\rho')$ is the union $\cup
\widetilde \Pi$ of $|\q\times \q'|$ parallelograms $\widetilde\Pi$
(shown doubly hatched in Figure~4) each of area $\rho \rho'/|\q \times
\q'|$.  By similarity, the ratio
$$
|\cup \tilde \Pi|:|\cup \Pi| = |\widetilde \Pi|:|\Pi| = 4\rho \rho':1.
$$
 But $|\cup \Pi| = 1$ and $|\cup \widetilde \Pi| 
= | B(\q;\rho) \cap B(\q',\rho')|$, whence
\begin{equation}
  \label{eq:Bli}
 | B(\q;\rho) \cap B(\q',\rho')|=4\rho \rho' =|B(\q;\rho)||B(\q';\rho')| 
\end{equation}
and $B(\q;\rho)$, $B(\q',\rho')$ are independent (see figure 4).
\begin{center}  
\begin{pspicture}(-7,-1)(6,6)
\psset{unit=0.5cm}
\psline(-5,10)(5,10) 
\psline(-5,0)(5,0) 
 \psline(-5,0)(-5,10) 
 \rput(-5,-0.5){$(0,0)$}
\rput(5,-0.5){$(1,0)$}
 \rput(-5,10.5){$(0,1)$}
 \rput(5,10.5){$(1,1)$}
  
\psline(5,0)(5,10) 

\psline(-5,9)(-3,10)  
\pspolygon[fillstyle=vlines](-5,9)(-5,10)(-3,10)

\psline(-5,6)(3,10) 
 \psline(-5,5)(5,10) 
 \psline(-5,4)(5,9) 
\pspolygon[fillstyle=vlines](-5,4)(5,9)(5,10)(3,10)(-5,6)

\psline(-5,1)(5,6) 
 \psline(-5,0)(5,5) 
 \psline(-3,0)(5,4) 
\pspolygon[fillstyle=vlines](-5,0)(-5,1)(5,6)(5,4)(-3,0)

\psline(3,0)(5,1)  
\pspolygon[fillstyle=vlines](3,0)(5,1)(5,0)

\rput(-0.25,5){$\Pi$}

\psline(4.5,10)(5,9.8)  
\pspolygon[fillstyle=hlines](4.5,10)(5,10)(5,9.8)

 \psline(5,.2)(-5,4.2) 
 \psline(5,0)(-5,4) 
 \psline(4.5,0)(-5,3.8) 

\pspolygon[fillstyle=hlines](5,0)(5,.2)(-5,4.2)(-5,3.8)(4.5,0)

 \psline(5,4.2)(-5,8.2) 
 \psline(5,4)(-5,8) 
  \psline(5,3.8)(-5,7.8) 

\pspolygon[fillstyle=hlines](5,4.2)(-5,8.2)(-5,7.8)(5,3.8)

 \psline(5,8.2)(0.5,10) 
  \psline(5,8)(0,10) 
  \psline(5,7.8)(-.4,10) 

\pspolygon[fillstyle=hlines](5,8.2)(0.5,10)(-.4,10)(5,7.8)

 \psline(.4,0)(-5,2.2) 
  \psline(0,0)(-5,2) 
  \psline(-.4,0)(-5,1.8) 

\pspolygon[fillstyle=hlines](.4,0)(-5,2.2)(-5,2.4)(-5,1.8)(-.4,0)

\psline(5,2.2)(-5,6.2)  
\psline(5,2)(-5,6)  
\psline(5,1.8)(-5,5.8) 

\pspolygon[fillstyle=hlines](-5,6.2)(5,2.2)(5,1.8)(-5,5.8) 

\psline(5,6.2)(-4.5,10)  
\psline(5,6)(-5,10)  
\psline(5,5.8)(-5,9.8) 

\pspolygon[fillstyle=hlines](-4.5,10)(5,6.2)(5,5.8)(-5,9.8)(-5,10)

 \psline(-4.5,0)(-5,.2) 

\pspolygon[fillstyle=hlines](-5,0)(-4.5,0)(-5,0.2)

\end{pspicture}

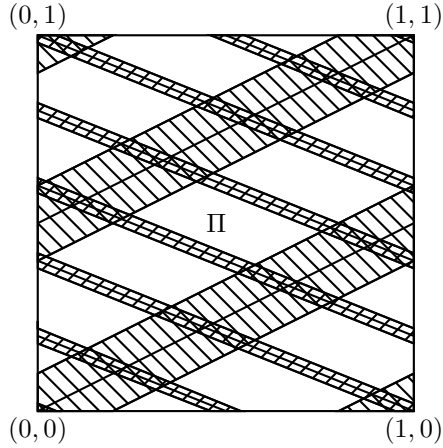
  \small{ Figure 4. $B(1,-2;\rho)\cap B(2,5;\rho')$ in 
    $\mathbb{T}^2$}
\end{center}

The extension to the general case again follows by considering $n $ copies
of the 2-dimensional space spanned by $\q$ and $\q'$ and the volume of
the corresponding $mn$-dimensional prisms. This gives 
for $m\ge2$, 
\begin{equation*}
  \label{eq:Bindep}
|\Bqbr\cap \Bqbrp|  = |\Bqbr|.|\Bqbrp|
 = 2^{2n}\rho^n {\rho'}^n. 
\end{equation*}
Thus when $m\ge2$ the pairwise probabilistically independent vectors
$\q$ in $\Z^n$ are precisely the pairwise linearly independent integer
vectors in $\Z^n$.  Thus to apply Theorem~\ref{thm:pwi2}, we need a
`large' set of such vectors and to find one, some number theoretic
ideas are required.

\subsection{A set of pairwise linearly independent vectors} 
The following argument is drawn from~\cite{Sprindzuk}. Let the highest
common factor of the integer components $q_1,\dots, q_m$ of
$\q\in\Z^m$ be denoted by $(\q)$, so that$(\q)=hcf(q_1,\dots, q_m)$.
The vector $\q\in\Z^m$ is said to be {\em primitive} if $(\q)=\pm
1$. Two distinct primitive vectors, $\q,\q'$ say, are linearly
independent when $\q,\q'$ are primitive. For if $\q,\q'$ are linearly
dependent, then $a\q = a'\q' $ for some real $a, a'$, $a,a'$ can be
assumed to be coprime integers (\ie, integers with no common factors
other than $\pm 1$).  Thus $a'$ divides each component $q_1,\dots,q_m$
of $\q$ and $a$ divides each component $q'_1,\dots,q'_m$ of $\q'$.
Since $\q,\q'$ are primitive, $a, a' = \pm 1$.

If in addition, $q_m,q'_m \ge 1$, then $a = a' = 1$ and $\q = \q'$. In
other words, no pair of distinct integer vectors in the set
$$
\cP_N = \{\q \in \Z^m:\q {\text{ primitive}}, |\q| \le N, q_m \ge 1\}
 = \{\q \in \Z^m:(\q)=1, |\q| \le N, q_m \ge 1\}
$$
is linearly dependent.  This set is the union of disjoint subsets (or
`hemispheres') $S_k$, consisting of vectors $\q$ in $\cP_N$ with
`radius' $|\q| = k$, i.e.,
$$
\cP_N = \bigcup_{k=1}^N S_k.
$$
Now let 
$$
\cP_{\infty} = \{\q \in \Z^m\col (\q)=1, q_m \ge 1\} = 
\bigcup_{k=1}^{\infty} S_k.
$$
Then distinct vectors $\q, \q' \in \cP_{\infty}$ are linearly
independent and so  $B_\q, B_{\q'}$ are independent, \ie,
$|B_\q\cap B_{\q'}|=|B_\q| |B_{\q'}|$.

The number of vectors $\q$ in $\Z^{m}$ with $|\q| = k$ and $q_m \ge 1$
is $2(m-1)(2k+1)^{m-2}k$, since each coordinate $q_j$, $j =
1,\dots,m$, satisfies $|q_j| \le k$ and $|q_{j'}| = k$ for some $j'$,
$1\le j' \le k$.  To obtain an asymptotic formula for $\Card S_k$,
divide up these vectors $\q$ in $\Z^m$ into classes $S(h)$ where $
h|k$ ($h$ divides $k$).  Then a vector $\q \in S(h)$ is of the form
$$
\q = (q_1,\dots,q_m) = (hr_1,\dots,hr_m) = h\rb,
$$
where $hr_{j'}  = q_{j'} = k$ and $\rb$ is primitive ($(\rb)=1$).  
Now
 \begin{align*}
   \Card S_k &  = \sum_{\q\in \cP_N, |\q|=k} 1 
= \sum_{(\q)=1, |\q| = k, q_m\ge 1} 1 
= \sum_{\substack{\q\in \Z^m\\ |\q|=k, q_m\ge1}} \sum_{d|(\q)}
   \mu(d) \\
&   =  \sum_{d|k}\mu(d) \sum_{\substack{|\rb|=k/d\in \cP_N\\
       r_m \ge 1}} 1 
 = \sum_{d|k}\mu(d)
   \left(2\frac{k}{d}+1\right)^{m-2}\frac{k}{d}2(m-1) \\ & =
   2^{m-1}(m-1) \sum_{d|k} \mu(d)(k/d)^{m-1} + O\left(k^{m-2}
     \sum_{d|k}\frac{|\mu(d)|}{d^{m-2}} \right),
\end{align*}
where  $\mu(d)$ is the M\"obius function~\cite[p 234]{HW}, given by
$$
\mu(d) = \begin{cases} (-1)^r &\text{when $d$ is the product of $r$ distinct 
primes,}\\
                 0     &\text{otherwise}, 
\end{cases}
$$
and has the important property that $\sum_{d|k}\mu(d)=1$ when $k=1$
and $0$ otherwise.  
 But $\vphi(k)$, the number of integers less than $k$ and coprime to
 $k$, is given by
$$
\vphi(k) := \sum_{1\le j\le k, (j,k)=1}1 =k \sum_{d|k} {\frac{\mu(d)}{d}},
$$
whence for $m=2$, 
\begin{equation*}
\Card S_k = 2 k\sum_{d|k} \mu(d)(1/d)
+ O\left(\sum_{d|k}|\mu(d)| \right) = 2 \vphi(k) + O(d(k))\asymp \vphi(k),
\end{equation*}
since $d(k)= \sum_{d|k} 1$, the number of divisors of $k$, satisfies
$d(k) = O(k^{\delta}) $ for any positive $\delta$~\cite[Theorem
315]{HW}.
 When the real part of the complex number $z>1$, Riemann's zeta
function is given by
$$
\zeta(z) = \sum_{k=1}^{\infty} \frac1{k^z}
 = \prod_{p{\text{ prime}}}\left(1-\frac1{p^{z}}\right)^{-1},
$$ 
so that for $m\ge 3$,
\begin{equation*}
  \frac1{\zeta(m-1)} =
  \left(\sum_{k=1}^{\infty} \frac1{k^{m-1}}\right)^{-1} 
< \prod_{{\substack{p|k\\ p\ \text{prime}}}}
  \left(1-\frac1{p^{m-1}}\right) = \sum_{d|k} \frac{\mu(d)}{d^{m-1}} <  1.
\end{equation*}
Thus
\begin{equation*}
\Card S_k \asymp \begin{cases} \vphi(k),   & m=2 \\
                     k^{m-1}, &  m \ge 3.
\end{cases}  
 \end{equation*}
 Now, as is well known, $\vphi(k)$ is comparable `on average' to
 $k$~\cite[Theorem 330]{HW} or to be precise,
\begin{equation}
  \label{eq:Phi}
  \Phi(N) = \sum_{k=1}^{N}\vphi(k) 
= \frac3{\pi^2} N^2 + O(N \log N) \asymp N^2 
\end{equation}
(see~\cite{HW,Sprindzuk}) 
and it turns out that if $\psi$ is non-increasing, the volume sum
determines the Lebesgue measure of $W(\psi)$.

As has been said, Fourier analysis of the periodic function
$\chi_{\Bqbr}$ can also be used; the linear independence of the $\q$'s
is crucial in establishing the measure of the intersection
$B(\q;\rho)\cap B(\q';\rho')$ (see~\cite[Chapter~1,\S5]{Sprindzuk} for
details).

\subsection{ Completing the proof of \K-Groshev theorem}

Next the divergence of the volume sum~\eqref{eq:volsum} is shown to
imply the divergence of a related sum over independent integer
vectors.
For each $N = 1,2,\dots,$
the partial sum
\begin{align*}
\sum_{\q \in \cP_{N}} |B_\q|
& = 2^n \sum_{\q \in \cP_{N}} \psi(|\q|)^n 
 = 2^n\sum_{k=1}^{N} \sum_{\q \in S_k} \psi(|\q|)^n 
 = 2^n \sum_{k=1}^{N} \psi(k)^n \sum_{\q \in S_k} 1 
 = 2^n \sum_{k=1}^{N} \psi(k)^n \Card S_k \\
& \asymp \begin{cases} \sum_{k=1}^N \vphi(k)\psi(k)^n & \text{when} \ m=2 \\
\sum_{k=1}^N k^{m-1}\psi(k)^n            & \text{when} \ m \ge 3.
\end{cases}
\end{align*}
Thus when $m\ge3$, the divergence of the sum~\eqref{eq:volsum} implies
the divergence of $\sum_{\q\in \cP_\infty} |B_\q|$.  To deal with the
case $m=2$, we use $\Phi(N) \asymp N^2 $~\cite{HW} and that the
monotonicity of $\psi(k)$ 
implies that $\sum_{k=1}^N \vphi(k)\psi(k)^n$ is comparable to
$\sum_{k=1}^N k\psi(k)^n$.
Hence if $\psi(k)$ is decreasing, then 
$$
\sum_{\q\in \cP_N} |B_\q| \asymp \sum_{k=1}^{N}  \vphi(k) \psi(k)^n   \asymp
\sum_{k=1}^{N} k \psi(k)^n 
$$
and the divergence of the right hand sum implies the divergence of the
left hand sum, which in turn implies $|W(\psi)|=1$. 

There are two interesting refinements of the \KG{} theorem.  The first
arises in the divergent case and is a quantitative version in the
sense of an asymptotic formula for the number of
solutions~\cite{HarmanMNT,Sprindzuk}. In the case of the real numbers,
the number $\cN(N;\alpha)$ of solutions with $q\le N$ of the
inequality
\begin{equation*}
    \left|\al-\frac{p}{q}\right| < \frac{\psi(q)}{q},
\end{equation*}
where again $\psi$ is decreasing, is
\begin{equation*}
  \cN(N;\alpha)=
2\sum_{q=1}^N \psi(q)(1 + o(1)).
\end{equation*}
For simultaneous \DA, where monotonicity can be omitted for dimensions
at least 2, an asymptotic formula holds for dimension at least
3~\cite{Gallagher65}. Asymptotic formulae will
not be discussed here but details are in~\cite{HarmanMNT,Sprindzuk}.
The other refinement concerns the finer structure of the null set when
the series converges.

\section{Hausdorff dimension}

It is a familiar fact that the one-dimensional Lebesgue measure of a
unit segment is 1 but that the two dimensional or planar Lebesgue
measure is 0.  This simple example illustrates that a dimension is
associated with the determination of Lebesgue measure.  For standard
shapes, such as the real line, the plane, rectangle or circle, the
dimension is the topological dimension and so is integral and obvious
(from our point of view, it could be called the Lebesgue dimension).
Of course, it is the Lebesgue measure of a set which is of
interest. As is well known, exceptional sets of $n$-dimensional \LM{}
zero can be studied using the more delicate notions of Hausdorff
dimension and measure, which allow null sets to be distinguished.
Hausdorff dimension, which is defined in terms of \HM, is a
generalisation of Lebesgue dimension and the two notions coincide for
standard sets. However, they differ in that any set in finite
dimensional Euclidean space has a Hausdorff dimension (which in
general will not be an integer).  In particular null sets have a \HD,
thus offering a way of studying sets that are `invisible' or
`negligble' in terms of Lebesgue measure and of distinguishing between
them. By contrast, a set of positive Lebesgue measure has full \HD{}
(equal to the Lebesgue dimension of the ambient space).

Although conceptually a simple but profound extension of
Carath\'eodory's construction of Lebesgue measure, Hausdorff measure has
a somewhat complicated definition and the reader is referred
to~\cite{MDAMshort,FalcGFSshort,FalcFG,Fed,MattilaGS,Rogers} for
fuller accounts.  For completeness a simpler formulation suited to our
purposes will be sketched.  Let $\mathcal C$ be a finite or countable
collection of open hypercubes $C \subset \R^k$ with sides of length
$\ell(C)$ and parallel to the axes.  For each non-negative real number
$s$ the $s$-volume of the collection $\cC$ is defined to be
$$
\ell^s(\mathcal C) = \sum_{C\in \mathcal C} \ell(C)^s.
$$ 
For any set $E$ in $\R^k$ and any real number $\delta>0$, let
$\cH_{\delta} ^s(E) = \inf \, \ell^s(\mathcal C_{\delta})$ be the
infimum taken over all `approximating' covers $\mathcal C_{\delta}$ of
$E$ by hypercubes $C$ with side length $\ell(C)$ at most $\delta$.  When $0
< \delta < 1$ and $t < s$, $\ell(C)^s < \ell(C)^t$ and so the number
$\cH_{\delta} ^s(E)$ decreases as $s$ increases.  The $s$-dimensional
outer measure $\cH^s(E)$ of $E$ defined by
$$
\cH^s(E) = \sup \{\cH_{\delta}^s(E) \col \delta > 0\} =\lim_{\delta\to
  0} \cH_{\delta}^s(E)
$$
is comparable to Hausdorff outer measure.  
If $t>s$, then $\ell(C)^t \le \delta^{t-s} \ell(C)^s$,   
whence
$\cH_{\delta}(E) \le \delta^{t-s} \cH_{\delta}^s(E)$, so that when 
$\cH^t(E)$ is positive, $\cH^s(E)$ is infinite and when $\cH^s(E)$
is finite, $\cH^t(E)$ vanishes.  The Hausdorff dimension $\dim E$ of $E$ is 
defined by
$$
\dim E = \inf \{s \in \R \col \cH^s(E) = 0 \},
$$
so that
\begin{equation*}
\cH^s(E) =   \begin{cases}  \infty, & s < \dim E, \\
                   0     & s > \dim E. 
  \end{cases}
  \end{equation*}
  Thus the dimension is that value of $s$ at which $\cH^s(E)$ `drops'
  discontinuously from infinity (see~figure 5).  Determining the \HM{}
  at this value is not alway easy and will not be discussed (but
  see~\S\ref{sec:FD} below). A cover for $E$ serves as a cover for any subset
  $E'$ of $E$ and so $E' \subseteq E$ implies that
$$
\dim E' \le \dim E.
$$

A comparison with viewing an object under a microscope can be made.
If the microscope lens is too close to the object $E$, the image fills
the eyepiece and cannot be resolved; if the lens is too far away, the
image is invisible.  At the focal length (\ie, at $s=\dim A$), the
image is in focus and can be seen properly. Thus the \HD{} is like
the focal length -- a \HM{} can be assigned to the set under
consideration at the \HD{} (see figure~5).
\begin{center}
  \begin{pspicture}(-4,-2)(6,6) 
\psline(-3,0)(4,0) 
  \rput(-3.2,-0.25){$(0,0)$}
  \rput(1.50,-0.25){$s$}
  \rput(-3.5,5){$\infty$}
\rput(-.250,-0.25){$\dim E$}
\rput(-3.75,3){$\cH^s(E)$}
 \pscircle[fillstyle=solid,fillcolor=black](0,2){.05} 
 \psline(-3.05,2)(-2.95,2)
 \rput(-3.9,2){$\cH^{\dim E}(E)$}
 \psline(-3,0)(-3,3.5) 
\psline(0,0)(4,0) 
 \psline(-3,3.6)(-3,3.8) 
 \psline(-3,3.9)(-3,4.04) 
 \psline(-3,4.15)(-3,4.25) 
\psline(-3,4.35)(-3,4.42) 
\psline(-3,4.5)(-3,4.55) 

\psline(-3,5)(0,5) 
\psline(0,0)(4,0) 

  \pscircle[fillstyle=solid,fillcolor=black](-2,-1){.5}
  \pscircle(0,-1){.5}
\pscircle[fillstyle=solid,fillcolor=gray](0.1,-.95){.25}
  \pscircle(2,-1){.5}
\end{pspicture}
\end{center}

\begin{center}
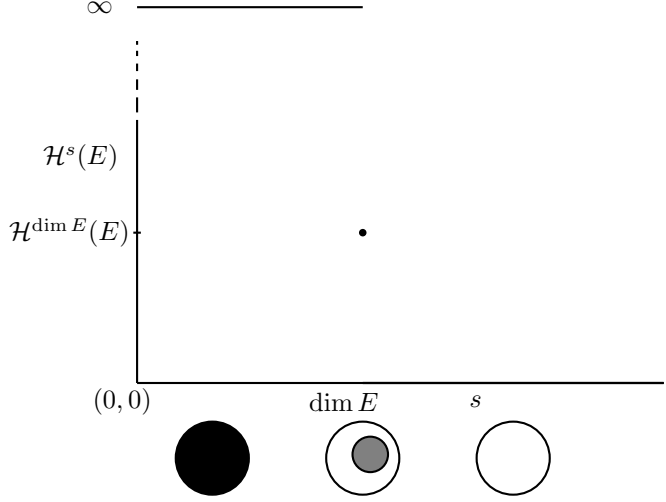
  \small{Figure 5. Graph of Hausdorff measure $\cH^s(E)$ against
    exponent $s$. The three circles represent \\
views through a microscope.}
\end{center}

When $s$ is a non-negative integer $m$ say, Hausdorff's $m$-measure is
comparable with Lebesgue's $m$-dimensional measure (and they agree
when $m=1$).  The Hausdorff dimension of a set $A \subseteq \R^k$ of
Lebesgue measure 0 is often established by obtaining an upper and a
lower inequality which combine to give the desired equality. In the
case of lim sup sets, such as $W(\psi)$, the upper bound usually
follows straightforwardly from a natural cover arising from the
definition and is closely related to the cover used in determining the
Lebesgue measure. For simplicity, let us take $\psi(r)=r^{-v}$, where
$v>0$, and write $W(\psi)=W_v$. Then for $W_v$, the cover $\cC$ of
hypercubes arises from the sets $B_{\q}$, where $ \q$ is a non-zero
integer vector and a straightforward calculation of the $s$-volume
$\ell^s(\cC)$ of the cover (see for example~\cite{mmd92}) gives
  \begin{equation}
    \label{eq:ub}
    \dim W_v\le
    \begin{cases}
      (m-1)n + \frac{m+n}{v +1} & \text{ when }
      v > \frac{m}{n} \\
      mn & \text{ when } v \le \frac{m}{n}.
    \end{cases}
  \end{equation}

  The lower bound is usually more difficult.  In the case of $W_v$,
  the argument can be shortened considerably by using the idea of
  `ubiquity'.  This was introduced originally to systematise and
  extend the determination of the lower bound for the Hausdorff
  dimension of sets of number theoretic and physical
  interest~\cite{DRV90a}.

\subsection{Ubiquity}
We start with some definitions and then introduce a lim sup set which
is associated with $W(\psi)$ and easier to work with.  Denote the
$\delta$-neighbourhood of the resonant set $\Rq$ by
$$
\widetilde{B}(\q;\delta) = \{X\in \T^{mn} \col \dist_\infty(X,\Rq)<\delta \}
$$
where $\dist_\infty(X,\Rq)= \inf\{\dist_\infty(X,U)\col U\in \Rq\}$ is
the distance in the supremum norm from $X$ to $\Rq$. This is not the
same set as the set 
\begin{equation*}
  B(\q;\delta):=\{X\in \T^{mn}\colon  \|\q X\|<\delta\},
\end{equation*}
which when $m=1$ reduces to the set $B(\q;\delta)$ defined by 
$\{\ub\in [0,1]^m\colon \|\q\ip \ub\| < \delta\}$. 
However it is readily shown that the sets are related by the following
inclusions: when $\q \ne {\mathbf0}$,
\begin{equation}
  \label{eq:inclusions}
  \widetilde{B}(\q;\frac{\delta}{m|\q|}) \subseteq \Bqbd 
  \subseteq \widetilde{B}(\q;\frac{\delta}{|\q|})   
\end{equation}
Let
$\tilde{\rho}\col\N \to (0,\infty)$ be a decreasing function.  When for the
family $\cR = \{R_{\q}:\q \in \Z^m \setminus \{{\mathbf0}\}\}$,
\begin{equation*}
|\T^{mn} \setminus \bigcup_{1\le|\q|\le N} 
\widetilde{B}(\q;\tilde\rho (N))| \to 0 \, {\text {as}}\, N \to \infty,
\end{equation*}
we say that $\cR$ is {\it{ubiquitous with respect to the function
    $\tr$}}.  

Let 
$$
\tr(N) = 2N^{-1-\frac{m}{n}} \log N. 
$$
When $m\ge 2$, the independence of the sets $B(\q;\rho)$ implies that
$\cR$ is ubiquitous with respect to $\tr$~\cite{GPIUsp} and a general
form of Dirichlet's theorem implies ubiquity without restriction on
the dimension $m$~\cite{mmd92,DRV90a} (see~\cite{HDSV97} for Hausdorff
measure results).  In essence this means that `most' $X$ are within
$\tr(N) = 2 N^{-1-m/n} \log N$ of some resonant set $R_{\q}$ with
$1\le|\q| \le N$.

Consider the lim sup set
\begin{equation*}
  \Lambda(\tp)=\left\{X\in \T^{mn}\colon \left|X-\Rq\right|
    < \tp(|\q|) \ \text{for infinitely many } \q\in\Z^m\right\} 
= \bigcap_{N=1}^{\infty} \bigcup_{|\q|=N}^{\infty}
  \widetilde{B}(\q;\tp(N)),
\end{equation*}
where $\tp(N)= N^{-v-1}/m$. 
By~\cite{mmd92,DRV90a}, the ubiquity of the family $\cR$
with respect to $\tr$ implies that the \HD{} of $\Lambda(\tp)$ 
satisfies
\begin{equation*}
  \dim \Lambda(\tp)\ge \dim \cR + \gamma \codim \cR,
\end{equation*}
where $\dim \cR$ is the topological dimension $(m-1)n$ of the resonant
set $\Rq$ and codimension $\codim \cR=n$ and
\begin{align*}
  \gamma &
  = \min\biggl\{1,\limsup_{N\to\infty}\frac{\log \tr(N)}{\log \tp(N)}\biggr\} \\
  &= \min\biggl\{ 1, (1+\frac{m}{n}) \frac1
  {\liminf_{N\to\infty}\frac{\log (mN^{1+v})}{\log N}}
  \biggr\} \\
  &= \min\biggl\{1, \frac{1+\frac{m}{n}}{1+v} \biggr\}.
  \end{align*}
Thus 
$$
\dim \Lambda(\tp)  \ge \min\bigl\{mn,\ (m-1)n + \frac{m+n}{v+1} \bigr\}. 
$$

By~\eqref{eq:inclusions} and the choice of $\tp$,
$$
\tB(\q;\tp(N)) =\tB(\q;N^{-v-1}/m) \subseteq B(\q;N^{-v}),
$$
whence
$$
\Lambda(\tp) \subset W_v. 
$$
Combining this with~\eqref{eq:ub} yields
\begin{equation*}
    \dim W_v=
    \begin{cases}
      (m-1)n + \frac{m+n}{v +1} & \text{ when }
      v > \frac{m}{n} \\
      mn & \text{ when } v \le \frac{m}{n}.
    \end{cases}
\end{equation*}
In one dimension, ubiquity is essentially equivalent to the `regular
systems' introduced by Baker and Schmidt~\cite{BS} and the above
result reduces to the \JB{} theorem.

\section{Further developments}
\label{sec:FD}
Determining the Hausdorff dimension of a set can be difficult enough
and finding the Hausdorff measure can be even harder without special
arguments available (such as when the Hausdorff measure coincides with
Lebesgue measure).  In another of his pioneering papers~\cite{Ja29},
\Ja{} established the \HM{} analogue of \K's theorem for simultaneous
\DA {} and showed that the Hausdorff $s$-measure at the critical
exponent (where $s=\dim W(\psi)$) is infinite. Dickinson and Velani
extended this result to systems of linear forms in~\cite{HDSV97}.
More recently with Beresnevich, they have developed a powerful and
unifying framework for obtaining the Hausdorff measure of lim sup sets
in the general setting of a compact metric space endowed with a
non-atomic probability measure and containing a family of resonant
sets~\cite{BDV06}.  The lim sup sets consist of points which lie close
to infinitely many resonant sets and include a very wide range of
results in the theory of metric \DA, including the set $W(\psi)$
discussed above. For recent applications, see the paper by
Drutu~\cite{drutu05} and the paper of Beresnevich and
Velani~\cite{BVParis07} in this proceedings. A similarity between the two
main theorems in~\cite{BDV06} suggests an equivalence between certain
Lebesgue and Hausdorff measure results and a Hausdorff measure
analogue of the Duffin-Schaeffer
conjecture~\cite{DuffinSchaeffer,PV90}. This is treated in a
subsequent paper by Beresnevich and Velani~\cite{BV06}.

\section{Acknowledgements}
I am grateful to Francoise Dalbo and Cornelia Drutu for organising the
most enjoyable conference on ``Dynamical systems and Diophantine
Approximation'', held at the Institut Henri Poincar\'e, Paris, 7-9
June 2004 and for an invitation to the Painlev\'e Institut,
Universit\'e de Lille 10-12 June. Their hospitality and that of the
institutes was much appreciated.  I am also grateful to Victor
Beresnevitch for his very helpful suggestions.

\bibliographystyle{amsplain}

\providecommand{\bysame}{\leavevmode\hbox to3em{\hrulefill}\thinspace}
\providecommand{\MR}{\relax\ifhmode\unskip\space\fi MR }
\providecommand{\MRhref}[2]{
  \href{http://www.ams.org/mathscinet-getitem?mr=#1}{#2}
}
\providecommand{\href}[2]{#2}

\end{document}